\documentclass[11pt]{article}
\usepackage{amssymb,amsfonts,latexsym}
\newcommand{\cleqn}{\setcounter{equation}{0}}
\newcommand{\clth}{\setcounter{theorem}{0}}
\newcommand {\sectionnew}[1]{\section{#1}\cleqn\clth}
\newcommand{\beq}{\begin{equation}}
\newcommand{\eeq}{\end{equation}}
\newcommand{\beqa}{\begin{eqnarray}}
\newcommand{\eeqa}{\end{eqnarray}}
\newcommand{\beaa}{\begin{eqnarray*}}
\newcommand{\eaa}{\end{eqnarray*}}
\newcommand{\nn}{\hfill\nonumber}
\newcommand{\text}{\textrm}
\newcommand \nc {\newcommand}
\nc \proof {{\em{Proof.\/}} }
\nc \qed {$\Box$\hfill}
\newtheorem{theorem}{Theorem}[section]
\newtheorem{lemma}[theorem]{Lemma}
\newtheorem{definition-theorem}[theorem]{Definition-Theorem}
\newtheorem{proposition}[theorem]{Proposition}
\newtheorem{corollary}[theorem]{Corollary}
\newtheorem{definition}[theorem]{Definition}
\newtheorem{example}[theorem]{Example}
\newtheorem{remark}[theorem]{Remark}
\nc \bth[1] { \begin{theorem}\label{t#1} }
\nc \ble[1] { \begin{lemma}\label{l#1} }
\nc \bdeth[1] { \begin{definition-theorem}\label{dt#1} }
\nc \bpr[1] { \begin{proposition}\label{p#1} }
\nc \bco[1] { \begin{corollary}\label{c#1} }
\nc \bde[1] { \begin{definition}\label{d#1}\rm }
\nc \bex[1] { \begin{example}\label{e#1}\rm }
\nc \bre[1] { \begin{remark}\label{r#1}\rm }
\nc \bcon[1] { \medskip\noindent{\it{Conjecture #1}} }
\nc \bqu[1]  { \medskip\noindent{\it{Question #1}} }
\renewcommand {\eth} { \end{theorem} }
\nc {\ele} { \end{lemma} }
\nc {\edeth}{ \end{definition-theorem} }
\nc {\epr} { \end{proposition} }
\nc {\eco} { \end{corollary} }
\nc {\ede} { \end{definition} }
\nc {\eex} { \end{example} }
\nc {\ere} { \end{remark} }
\nc {\econ} {\smallskip}
\nc {\equ} {\smallskip}
\nc \eqref[1] {{\rm{(\ref{#1})}}}
\nc \thref[1]{Theorem \ref{t#1}}
\nc \leref[1]{Lemma \ref{l#1}}
\nc \prref[1]{Proposition \ref{p#1}}
\nc \coref[1]{Corollary \ref{c#1}}
\nc \deref[1]{Definition \ref{d#1}}
\nc \exref[1]{Example \ref{e#1}}
\nc \reref[1]{Remark \ref{r#1}}
\def \lg {\langle}
\def \rg {\rangle}
\def \g  {\mathfrak{g}}
\def \f  {\mathfrak{f}}
\def \x  {\xi}
\def \e  {\eta}
\def \S  {S}
\def \s  {s}
\def \b  {\mathfrak{b}}
\def \a  {\mathfrak{a}}
\def \k  {\mathfrak{k}}
\def \n  {\mathfrak{n}}

\def \h  {\mathfrak{h}}
\def \l  {\mathfrak{l}}
\def \d  {\mathfrak{d}}
\def \mcA {\mathcal{A}}
\def \mcB {\mathcal{B}}
\def \wt  {\widetilde}
\def \tmcB {{\widetilde{\mathcal{B}}}}
\def \tB {\widetilde{B}}
\def \tW {{\widehat{W}}}
\def \D  {\mathcal{D}}
\def \o  {\otimes}

\def \op {\oplus}
\def \th {\theta}
\def \Th {\Theta}
\def \sub{\subset}
\def \sup{\supset}
\def \de  {\delta}
\def \sig {\sigma}
\def \ga {\gamma}
\def \al {\alpha}
\def \be {\beta}
\def \Om {\Omega}
\def \Rset {{\mathbb R}}
\def \Cset {{\mathbb C}}
\def \ra {\rightarrow}
\def \mt {\mapsto}
\def \ol {\overline}
\def \st {\ast}
\def \Lie { {\mathrm{Lie}} }
\def \id { {\mathrm{id}} }
\def \Ad { {\mathrm{Ad}} }
\def \ad { {\mathrm{ad}} }

\def \Ker { {\mathrm{Ker}} }
\renewcommand \Im { {\mathrm{Im}} }

\begin{document}
\title{The Double and Dual of a Quasitriangular Lie Bialgebra}
\author{Timothy J. Hodges \thanks{The first author was partially 
supported by grants from the National Security Agency and the Charles
P. Taft Foundation}\\Department of Mathematics
\\University of Cincinnati\\ Cincinnati, OH 45221-0025
\\timothy.hodges@uc.edu\\ \and Milen Yakimov 
\thanks{The second author was supported by NSF grants 
DMS94-00097 and DMS96-03239}
\\
Department of Mathematics \\
University of California at Berkeley\\
Berkeley, CA 94720\\
yakimov@math.berkeley.edu\\}

\date{}
\maketitle
\begin{abstract} Let $G$ be a connected, simply connected Poisson--Lie
group with quasitriangular Lie bialgebra $\g.$ An explicit description of
the double $\D(\g)$ is given, together with the embeddings of $\g$ and
$\g^\st$. This description is then used to provide a construction of the
double $\D(G).$ The aim of this work is to describe $\D(G)$ in sufficient
detail to be able to apply the procedures of  Semenov-Tian-Shansky and 
Drinfeld for the classification of symplectic leaves and Poisson
homogeneous spaces for Poisson--Lie groups.
\end{abstract}
\sectionnew{Introduction}
Knowing explicitly the Lie algebra structure of the double of a
Lie bialgebra is necessary for the classification
of the Poisson homogeneous spaces and symplectic leaves of the
corresponding Poisson--Lie groups.

To be more precise, let $\g$ be a Lie bialgebra over a field $k$
and
$\D(\g)$ be its double. By the dual Lie bialgebra $\g^\st$ of $\g,$ we will
mean the standard dual Lie bialgebra of $\g,$  equipped with opposite
cobracket. When $k=\Rset$ or $\Cset$
the connected, simply connected Poisson--Lie groups with tangent Lie
bialgebras
$\g,$ $\g^\st$ and $\D(\g)$ will be denoted by $G,$ $G^\st,$ and $\D(G).$
$\D(G)$ and $G^\st$ are called the double and dual Poisson--Lie groups of
$G,$ respectively.
Let $j \colon G \to \D(G)$ and $j^\st \colon G^\st \to \D(G)$ be the group
homomorphisms corresponding to the embeddings $\g \sub \D(\g)$ and
$\g^\st \sub \D(\g).$

According to a result of Drinfeld \cite{D}, the Poisson homogeneous spaces
for $G$ with connected stabilizers are in one to one correspondence with
the orbits of the adjoint action of $G$ on the variety of
the Lagrangian subalgebras $\l$ of $\D(\g)$ for which $\l \cap \g$
integrates to a closed subgroup of $G.$

By a Theorem of Semenov-Tian-Shansky \cite{STS2},
the symplectic leaves of $G$ are the orbits of the (local) dressing
action of $G^\st$ on $G$ which is the lifting
under $G \to \D(G)/j^\st(G^\st)$ of the left action of $G^\st$ on
$\D(G)/j^\st(G^\st).$ Again to work with it in a concrete situation, one
needs an explicit description of  $\D(G)$ and of the images of $G$ and
$G^\st$. For $G$ semisimple, the symplectic leaves of $G$ were described
using
this procedure by Hodges and Levasseur \cite{HL1} for the standard bialgebra
structure and by Yakimov \cite{Y} for the structures associated to a general
Belavin-Drinfeld $r$-matrix.

If the Lie bialgebra $\g$ is quasitriangular with an $r$-matrix
$r,$ the structure of $\D(\g)$ is encoded in the structure of the Lie
algebra $\g$ and the $r$-matrix. For a factorizable Lie bialgebra it was
shown by
Reshetikhin and Semenov-Tian-Shansky \cite{RS} that $\D(\g) \cong \g
\oplus
\g$.  Stolin \cite{Sto} and  Levendorskii, Soibelman \cite{LS} showed for
triangular structures on
complex simple Lie algebras and their compact real forms that
$\D(\g) \cong
\g\ltimes \g^\st$. 
These results enable one to construct $\D(G)$ as $G \times G$ and 
$G \ltimes \g^\st$ respectively.
The goal of this note is to extend these results to the general case; that
is,
to give  {\em{an explicit description of $\D(\g)$
and of the embeddings $\g, \, \g^\st \sub \D(\g)$ in terms of the
$r$-matrix
and then to lift this description to a construction of $\D(G)$.}}

Next we briefly explain our results. The subspace $\f$ of $\g$
spanned by all components of the tensor $r + r^{21} \in S^2(\g)$
is an ideal of $\g.$ In Section~3 we show that $\D(\g)$ is an extension of
the semidirect sum $\g \ltimes \f$ by an abelian Lie algebra.
The 2-cocycle $\al$ defining the extension is computed in Section~4.
Section~5 deals with the image of $\g^\st \sub \D(\g)$ under
this isomorphism. It is described in terms of an explicit 1-cochain
for a subalgebra $\b$ of $\g \ltimes \f$ whose coboundary is
essentially the restriction of $\al$ to $\b.$ In Section 6
we lift this construction to a construction of $\D(G)$ as an extension of
$G \ltimes F$ by $\f^\perp$, generalizing the above descriptions of
$\D(G)$
in the factorizable and triangular cases.
In Section ~7 we specialize the results to the two limiting cases of
triangular and factorizable structures, i.e. $r+r^{21}=0$ or of full rank. 

No assumptions will be made on the ground field $k$ and unless otherwised
specified all Lie bialgebras can be infinite dimensional. For a vector
space $V$ over $k,$ the natural pairing between $\x \in V^\st$
and $v \in V$ will be denoted both by $\x(v)$ and
$\lg \x, v \rg.$ The symbol $\ltimes$ will be used for both
semidirect sums and products of Lie algebras and Lie groups.

{\flushleft{\bf{Acknowledgements}} The authors would like to thank
Nicolai Reshetikhin and Sei-Qwon Oh for fruitful discussions.}

\medskip\noindent

\sectionnew{Lie bialgebras}
Recall that a Lie bialgebra $\g$ with Lie bracket $[.,.]$
and Lie cobracket $\de$ is called
quasitriangular
if there exists an element  $r \in \g \o \g$
such that
\[
\de(x) = (\ad_x \o \id + \id \o \ad_x) r, \, \forall x \in \g
\]
and $r$ satisfies the Yang--Baxter equation
\[
[r_{12}, r_{13}] + [r_{12}, r_{23}] + [r_{13}, r_{23}] =0.
\]

Let $\g^\st$ be the dual Lie bialgebra of $\g$ with opposite cobracket as
explained in the introduction. (This convention is not essential
since in Sections~3--7 we only work with the Lie algebra structure of
$\g,$ $\g^\st,$ and $\D(\g).)$
There are two natural Lie bialgebra homomorphisms
$r_\pm \colon \g^\st \to \g$ (see \cite{RS}), given by
\beqa
&& r_+(\x) = (\x \o \id) r,
\label{r+} \\
&& r_-(\x) = - (\id \o \x) r, \; \x \in \g^\st.
\label{r-}
\eeqa

Recall also that the double Lie bialgebra $\D(\g)$ of $\g$ is
characterized by the properties

1. $\D(\g) \cong \g \op \g^\st$ as a coalgebra,

2. $\g$ and $\g^\st$ are subalgebras of $\D(\g)$ and two elements
$x \in \g \sub \D(\g)$ and $\x \in \g^\st \sub \D(\g)$ commute
as
\beq
[x, \x] = \ad^\st_x(\x) - \ad^\st_\x(x)
\label{ad}
\eeq
where $\ad^\st$ denotes the coadjoint actions of both $\g$ and
$\g^\st.$

There is a canonical nondegenerate invariant bilinear form on $\D(\g)$
given by
\beq
\lg \lg x_1+ \x_1, x_2 + \x_2 \rg \rg :=
\lg \x_1, x_2 \rg +
\lg \x_2, x_1 \rg, \; \forall x_i \in \g, \, \x_i \in \g^\st
\label{in}
\eeq
with respect to which $(\D(\g), \g, \g^\st)$ is a Manin triple.
For details we refer to \cite{CP}.
\sectionnew{An exact sequence}
Let $\f \sub \g$ denote the span of all components of
the tensor
\beq
\Om = r + r_{21} \in S^2(\g).
\label{Om}
\eeq
Since $\de$ is skewsymmetric, $\Om$ is $\ad_{\g}$ invariant.
This implies that $\f$ is an ideal of $\g.$
Denote by $\f^\perp$ the subspace of $\g^\st$ orthogonal to $\f \sub \g:$
\beq
\f^\perp = \{ \e \in \g^\st \mid \lg \e , y \rg = 0, \forall \,
y \in \f \}.
\label{l'}
\eeq
In terms of the maps $r_\pm,$ $\f$ and $\f^\perp$ are given by
\beqa
&& \f = \Im (r_+ - r_-),
\label{charl} \\
&& \f^\perp = \Ker (r_+ - r_-).
\label{charl'}
\eeqa

Define two linear maps
\[
p_\pm \colon \D(\g) \to \g
\]
by
\beq
p_\pm (x + \x) = x + r_\pm(\x); \; x \in \g \sub \D(\g), \,
                                 \x \in \g^\st \sub \D(\g).
\label{ppm}
\eeq

 The following result is due to Majid \cite[8.2.7]{Mb}; we include a proof
for completeness.

\bpr{p_pm}The maps $p_\pm \colon \D(\g) \to \g$ are homomorphisms of Lie
algebras.
\epr
In fact the maps $p_\pm$ are homomorphisms of Lie bialgebras since
$\D(\g) \cong \g \op \g^\st$ as a coalgebra, but we will not need this
fact
here.

In the proof of this proposition and below we will use the following
simple fact.

\ble{coadj} For all $\x \in \g^\st$ and $x \in \g$
\beqa
\ad_\x^\st(x) &=&
-(\x \o \id) ( (\ad_x \o \id + \id \o \ad_x)r )
\nn \\
&=& (\id \o \x) ( (\ad_x \o \id + \id \o \ad_x)r ).
\nn
\eeqa
\ele
\proof For an arbitrary $\x_0 \in \g^\st$
\beqa
&& \lg \ad_\x^\st(x), \x_0 \rg =
-  \lg x, [\x, \x_0]       \rg
\nn \\
&& = - \lg \de(x), \x \o \x_0 \rg =
\lg -(\x \o \id) ((\ad_x \o \id + \id \o \ad_x)r ), \x_0 \rg.
\nn
\eeqa
\hfill \qed \\
{\em{Proof of \prref{p_pm}.}} The restrictions of $p_\pm$
to $\g \sub \D(\g)$ and $\g^\st \sub \D(\g)$
coincide with $r_\pm$ and thus are homomorphisms.
We are left with checking that
\[
[p_\pm (x), p_\pm(\x)] = p_\pm([x, \x])
\]
for $x \in \g \sub \D(\g),$ $\x \in \g^\st \sub \D(\g).$
Using \leref{coadj}, for the map $p_+$ one computes
\beqa
p_+([x, \x]) &=& p_+(\ad^\st_x(\x) - \ad^\st_\x(x) )
\nn \\
&=& (\ad^\st_x(\x) \o \id) r - \ad^\st_\x(x)
\nn \\
&=& - (\x \o \id) ((\ad_x \o \id) r)
    + (\x \o \id) ((\ad_x \o \id + \id \o \ad_x) r)
\nn \\
&=& [x, r_+(\x)] = [p_+(x), p_+(\x) ]
\nn
\eeqa
\hfill \qed

The kernels of $p_\pm$ are
\beq
\Ker \, p_\pm = (\id - r_\pm) \g^\st.
\label{kerd}
\eeq

\ble{commute} The algebras $\Ker \, p_+$ and $\Ker \, p_-$ mutually
commute in $\D(\g).$
\ele
\proof Since $\Ker \, p_\pm$ are ideals of $\D(\g),$
for any $\x_1, \x_2 \in \g^\st$
\[
[(\id - r_+)\x_1, (\id - r_-)\x_2] =
(\id - r_+) \x
\]
for some $\x \in \g^\st.$ Comparing the components in $\g^\st \sub \D(\g)$
of the two sides of this equality, we get
\beq
\x = [\x_1, \x_2] -
\ad^\st_{r_+(\x_1)}(\x_2) +\ad^\st_{r_-(\x_2)}(\x_1).
\label{r+-}
\eeq
Thus \leref{commute} is equivalent to the vanishing of the above
expression for all $\x_1, \x_2 \in \g^\st.$
For an arbitrary $x \in \g$
\beqa
&&\lg [\x_1, \x_2] - \ad^\st_{r_+(\x_1)}(\x_2) +\ad^\st_{r_-(\x_2)}(\x_1),
x \rg
\nn \\
&&=\lg \x_1 \o \x_2, \de(x) \rg + \lg \x_2, \ad_{r_+(\x_1)}(x) \rg -
\lg \x_1, \ad_{r_-(\x_2)}(x) \rg
\nn \\
&&=\lg \x_1 \o \x_2, \de(x)- (\id \o \ad_x)r - (\ad_x \o \id)r \rg = 0.
\nn
\eeqa
\hfill \qed

Combine the maps $p_\pm$ into the single morphism
\beq
p \colon \D(\g) \to \g \op \g, \; p(d)= (p_+(d), p_-(d)).
\label{p}
\eeq
Using \eqref{charl}, for its image we obtain
\beqa
\Im \, p &=& \{(x + r_+(\x), x + r_-(\x)) \mid x \in \g, \, \x \in
\g^\st
\}
\nn \\
      &=& \{(x_0, x_0 + y) \mid x_0 \in \g, \, y \in \f \}
\nn \\
      & \cong & \g \ltimes \f.
\label{gl}
\eeqa
In the semidirect sum
$\g \ltimes \f,$ $\g$ acts on $\f$
via the adjoint action (recall that $\f$ is an
ideal of $\g).$ The resulting map $\D(\g) \to \g \ltimes \f,$
given by $ x + \x \mt (x + r_+(\x), -(r_+ - r_-)\x),$
$x \in \g,$ $\x \in \g^\st$ will still be denoted by $p$. This map was also
constructed and discussed by Majid \cite[8.2.6]{Mb}.

Similarly, from \eqref{charl'} we see that the kernel of $p$ is
\beqa
\Ker \, p &=& \Ker \, p_+ \cap \Ker \, p_-
\nn \\
&=& (\id - r_+) ( \Ker (r_+ - r_-) )
\nn \\
&=& (\id - r_+) \f^\perp.
\eeqa
\leref{commute} implies that $\Ker \, p$ is an abelian subalgebra of
$\D(\g).$

Equip the space $\f^\perp$ with the structure of abelian Lie algebra.
Then the map
\beq
i= (\id - r_+)|_{\f^\perp} \colon \f^\perp \to \D(\g)
\label{i}
\eeq
is an embedding of Lie algebras.

To summarize:

\bpr{exact} The sequence
\beq
0 \to \f^\perp \stackrel{i}{\to} \D(\g) \stackrel{p}{\to}
\g \ltimes \f \to 0
\label{exact}
\eeq
is exact.
\epr

Note that $\f^\perp,$ considered as a subspace of $\g^\st,$ is
not an abelian Lie algebra in general. For example if $\g$
is a triangular Lie bialgebra, $\f^\perp = \g^\st$ but
$\g^\st$ is not always abelian.
\sectionnew{Description of the extension}
According to the previous section, $\D(\g)$ is an extension of
$\g \ltimes \f$ by the abelian Lie algebra $\f^\perp$. Thus, there is an
induced representation of $\g \ltimes \f$ on $\f^\perp$ and the extension is
uniquely defined by an element of $H^2(\g \ltimes \f, \f^\perp)$ (see, for
instance, \cite[XIV.5]{CE} or \cite[VII.4]{HS}). To describe these
structures, we need to fix a splitting 
$\S \colon \g \ltimes \f \to \D(\g)$  of \eqref{exact}
in the category of vector spaces, i.e. $p \circ \S = \id_{\g \ltimes \f},$
such that S(0)=0.
 {} From \eqref{gl} one sees that such splitting may be defined
by a right inverse map $\s \colon \f \to \g^\st$ of
$(r_+-r_-) \colon \g^\st \to \f$ (i.e. $(r_+ - r_-) \circ \s = \id_\f).$
The desired splitting $\S \colon \g \ltimes \f \to \D(\g)$ can then be
defined
by
\beq
\S( x,y) = x + r_+ \s(y) - \s(y) \in \D(\g), \, x \in \g, \, y \in \f
\label{def_S}
\eeq
(the elements of $\g \ltimes \f$ will be represented as ordered pairs).

\ble{one-rplus} The maps $\id-r_\pm \colon \g^\st \to \D(\g)$ are
$\g$-module
homomorphisms. That is, $(\id-r_\pm)(\ad_{x}^\st \x) = [x, (\id-r_\pm)\x]$
for
all $x \in \g$ and all $\x \in \g^\st$.
\ele

\noindent \proof
 In the $r_+$ case we have
\beaa[x, (\id-r_+)\x] &=& -[x, r_+\x] + \ad_{x}^\st \x - \ad_{\x}^\st x\\
 &=& -(\x \o \id) ( \id \o \ad_x)r +\ad_{x}^\st \x - \ad_{\x}^\st x\\
 &=& \ad_{x}^\st \x + (\x \o \id) ( \ad_x \o \id)r\\
 &=& (\id-r_\pm)(\ad_{x}^\st \x).
\eaa
The argument for $r_-$ is analogous. \hfill \qed

\ble{ext-action} The action of $\g \ltimes \f$ on $\f^\perp$ induced from
the extension
\[
0 \to \f^\perp \stackrel{i}{\to} \D(\g) \stackrel{p}{\to}
\g \ltimes \f \to 0
\]
is given by $(x, y)\cdot \e = \ad_{x}^\st(\e)$ for $x \in \g$, $ y \in \f,$
and $\e \in \f^\perp$.
\ele

\noindent\proof
 The action is given by
$(x, y)\cdot \e = i^{-1} ( [S(x,y), i(\e)])$. Now \leref{commute} implies
that $i(\f^\perp)$ is in the center of $\Ker \, p_+$ and therefore
\beaa
[\S( x, y ), i(\e)] &=& [x -(\id-r_+)s(y), (\id-r_+) \e] 
= [x,(\id-r_+)\e]\\
&=& (\id - r_+)( \ad_{x}^\st(\e)).
\eaa
by \leref{one-rplus}.
\hfill \qed

 Let us briefly recall the construction of an extension of a Lie algebra
$\h$ by an $\h$-module $V$ associated to an element of $H^2(\h, V)$. Let
$\alpha \in Z^2(\h, V)$ be a representative of the cohomology class
$[\alpha] \in H^2(\h, V)$. Then the extension of $\h$ by $V$ associated to
$\alpha$ is the vector space $\h \oplus V$ equipped with the bracket
\[
[(h_1, v_1),(h_2,v_2)] = ([h_1,h_2], 
  h_1 \cdot v_2 - h_2 \cdot v_1 + \alpha(h_1,h_2)).
\]
We will denote such an extension by $\h \ltimes_{\alpha}V$.

 Returning to $\D(\g)$, notice that it follows easily from
\leref{one-rplus}
that for any $x \in \g$ and $y \in \f$,
\[
\ad^\st_{x}(\s(y)) - \s([x, y]) \in \f^\perp.
\]
Hence we may define a 2-form, 
$\alpha \colon (\g \ltimes \f) \otimes (\g \ltimes \f) \to \f^\perp$ by
\beqa
\al((x_1, y_1), (x_2, y_2)) &=&
- \ad_{x_1}^\st (\s(y_2))
+ \ad_{x_2}^\st (\s(y_1))
- \ad_{y_1}^\st (\s(y_2))
\nn \\
&+& \s([x_1, y_2] - [x_2, y_1] + [y_1, y_2]).
\label{bedef}
\eeqa

\bth{extension} The form $\alpha$ is a 2-cocycle and
\beq
\label{ident}
 \D(\g) \cong (\g \ltimes \f)\ltimes_{\alpha} \f^\perp.
\eeq
\eth

\noindent\proof The cohomology class associated to the extension
\[
0 \to \f^\perp \stackrel{i}{\to} \D(\g) \stackrel{p}{\to}
\g \ltimes \f \to 0
\]
is the class of the 2-cocycle
\[
 i^{-1}([\S(x_1, y_1), \S(x_2, y_2)] -\S([(x_1, y_1),(x_2, y_2)])).
\]
But
\beaa
&& [\S(x_1, y_1), \S(x_2, y_2)] = [x_1 -(\id-r_+)s(y_1),
                                   x_2 -(\id-r_+)s(y_2)] 
\\  
&& = [x_1,x_2] + (\id - r_+)
   ( 
    - \ad^\st_{x_1}(\s(y_2))
    + \ad^\st_{x_2}(\s(y_1))- \ad^\st_{y_1}(\s(y_2))  )
\eaa
and
\[
 S([(x_1, y_1),(x_2, y_2)]) = [x_1,x_2] -(\id-r_+)\big([x_1,y_2]
 -[x_2,y_1] +[y_1,y_2]\big)
\]
So
\[
[\S(x_1, y_1), \S(x_2, y_2)] -\S([(x_1, y_1),(x_2, y_2)]) 
=(\id-r_+)(\al((x_1, y_1), (x_2, y_2)))
\]
as required.
 \hfill \qed

\bre{inv_form} In the case when $\g$ admits a nondegenerate invariant
bilinear form, $\g^\st$ can be identified as a linear space with $\g$ and 
$r_\pm$ can be thought of as linear endomorphisms of $\g.$  
In \thref{extension} $\f^\perp$ may be replaced by the 
orthogonal complement to $\f$ in $\g$ and the coadjoint representation
$\ad^\st$ of $\g$ with the adjoint one. The results of the
following section can be reformulated in a similar way.
\ere
\sectionnew{The embeddings of $\g$ and $\g^\st$ in $\D(\g)$}
It is clear that in the identification \eqref{ident},
$\g$ is embedded in $\D(\g)$ by
\beq
x \mt (x, 0) \ltimes_\al 0 \, \in (\g \ltimes \f) \ltimes_\al \f^\perp,
\; x \in \g
\label{embed-g}
\eeq
(see eqs. \eqref{i} and \eqref{def_S}).

 {}  From the definition of $\s$ it follows that for any $\x \in \g^\st$
\[
\x - \s (r_+ - r_-) \x \in \f^\perp.
\]
This easily gives
\[
\x = \S(r_+(\x), - (r_+ - r_-)\x) + i( \x - \s (r_+ - r_-)\x), \,
\forall \x \in \g^\st.
\]
Thus the embedding of $\g^\st$ in 
$(\g \ltimes \f) \ltimes_\alpha \f^\perp \cong \D(\g)$ is given by
\beq
\x \mt (r_+(\x), - (r_+ - r_-)\x) \ltimes_\al ( \x - \s (r_+ - r_-)\x).
\label{embed^st}
\eeq
In this section we describe in an invariant way the image of
this embedding.

It will be more convenient to work with the algebra
\beq
\a := \{ (x_1, x_2) \mid x_i \in \g,
 \, x_2 - x_1 \in \f \} \sub \g \op \g
\label{def_d}
\eeq
isomorphic to $\g \ltimes \f.$ The corresponding isomorphism as shown in
\eqref{gl}
is $(x_1, x_2) \mt (x_1, x_2 - x_1).$
In terms of $\a$ the identification \eqref{ident} is
\beq
\D(\g) \cong \a \ltimes_\alpha \f^\perp
\label{ident2}
\eeq
and the corresponding action of $\a$ on $\f^\perp$ is
\beq
(x_1, x_2) \cdot \e = \ad^\st_{x_1}(\e) \in \f^\perp, \;
(x_1, x_2) \in \a, \, \e \in \f^\perp.
\label{rep_d}
\eeq
It is also straightforward to write $\al$ as a 
2-cocycle for $\a$ with values in $\f^\perp.$ From \eqref{embed^st} it
follows that  $\g^\st$ is embedded in
$\a \ltimes_\alpha \f^\perp \cong \D(\g)$ by
\beq
\x \mt (r_+(\x), r_-(\x)) \ltimes_\al ( \x - \s (r_+ - r_-)\x), \, \x \in
\g^\st.
\label{embed2}
\eeq

We first state a general classification of subalgebras of extensions
of the form $\h \ltimes_{\alpha} V$ where $V$ is an $\h$-module and $\alpha
\in Z^2(\h, V)$. Let $p \colon \h \ltimes_{\alpha} V \to \h$ and 
$q \colon \h \ltimes_{\alpha} V \to V$ be the projections $p(h,v)=h$ and
$q(h,v)=v$. (Recall from the previous section that the map \eqref{p} is
a special case of the first projection.)

\ble{b-beta-W}
 Let $\b \subset \h$ and let $W$ be a $\b$-submodule of $V$. Let $\beta
\colon \b \to V/W$  be a 1-cochain whose coboundary is
$-\sigma \circ \alpha|_\b$ where $\sig$ denotes the
projection $V \to V/W.$ Define
\[
 \b^\beta_W = \{ (x, v) \mid x \in \b, \; v+W = \beta(x) \}.
\]
Then $\b^\beta_W$ is a Lie subalgebra of $\h \ltimes_{\alpha} V$.
\ele

\bth{subalgebras}
If $\k$ is a Lie subalgebra of $\h \ltimes_\alpha V$, then $\k$ is of the
form $\b^\beta_W$ where $\b = p(\k)$, $W=\k \cap V$ and $\beta \colon \k
\to V/W$ is given by $\beta(k) = q(p^{-1}(k)) + W$.
\eth

 We now identify the image of $\g^\st$ inside $\a \ltimes_\alpha \f^\perp$.
It follows from \eqref{embed^st} that 
$W = \Ker \, r_+ \cap \Ker \, r_-$. The
projection $\b =p(\g^\st)$ can be given the following conceptual
description.

Define $\g_\pm = \Im \, r_\pm$ and $\n_\pm = r_\pm( \Ker \, r_\mp).$
Recall from \cite{RS, STS} that $\n_\pm$ are ideals of
$\g_\pm.$
The map $\th \colon \g_+/\n_+ \to \g_-/\n_-$ given by
\[
\th( r_+(\x) + \n_+) = r_-(\x) + \n_-, \; \x \in \g^\st
\]
is well defined and is an isomorphism of Lie algebras.
It is called the Cayley transform of the
$r$-matrix $r.$ The projections from $\g_\pm$ to
$\g_\pm/\n_\pm$ will be denoted by $\pi_\pm.$

Let $\b$ be the subalgebra of $\a \cap (\g_+ \op \g_-) \sub \g \op \g$
consisting
of the elements $(x_+, x_-) \in \g_+ \op \g_-$ such that
\beq
\th \circ \pi_+(x_+) = \pi_- (x_-).
\label{def_a}
\eeq
The following lemma is a direct consequence from the definition of
$\th$ (cf. \cite{STS}).
\ble{p(gst)} The image of $\g^\st \sub \D(\g)$ under the projection
$p \colon \D(\g) \to \a$ coincides with $\b \sub \a:$
\[
p(\g^\st) = \b.
\]
\ele

 The projection $p|_{\g^\st} \colon \g^\st \to \a$ factors through to an
isomorphism $\overline{p} \colon \g^\st/W \to \b$ given by
\[
\ol{p}(\x + W) =(r_+(\x),r_-(\x)).
\]
Hence the 1-cochain $\beta \colon \b \to \ol{\f}^\perp$ of
\thref{subalgebras} in this situation will be given by
\beqa
 \beta(x_1,x_2) &=& \ol{p}^{-1}(x_1, x_2) - 
\sig s (r_+-r_-)\ol{p}^{-1}(x_1, x_2)\nn\\
 &=& \ol{p}^{-1}(x_1, x_2) - \sig s(x_1-x_2).
\label{betadef}
\eeqa

This description can be summarized as follows.

\bth{g^st}
 The image of $\g^\st$ inside $\D(\g)$ under the isomorphism $\D(\g)
\cong\a\ltimes_\alpha \f^\perp$ coincides with the subalgebra $\b^\beta_W$
where $\b$ is the algebra described in \leref{p(gst)}, $W= \Ker\; r_+ \cap
\Ker\; r_-$ and $\beta \colon \b \to \f^\perp/W$ is as described
in \eqref{betadef}.
\eth

\bre{Manintrip} To describe the Manin triple
$(\D(\g), \g, \g^\st)$ completely, we are left with
expressing the canonical
bilinear form \eqref{in} on $\D(\g)$ in terms of the identification
$\D(\g)\cong \a \ltimes_\al \f^\perp.$ This is needed for example for the
description of the Poisson homogeneous spaces for the Poisson--Lie groups
with tangent Lie bialgebra $\g$ as explained in the introduction.
For any
\[
d=(x_1, x_2) \ltimes_\al \e \in \a \ltimes_\al \f^\perp
\]
set
\beqa
&&\x(d) := \e + \s (r_+ - r_-) (x_1 - x_2),
\nn \\
&&x(d) := x_1 - r_+ (x(d)).
\nn
\eeqa
Using equations \eqref{embed-g} and \eqref{embed2}, it is easy to show that
the composition of the isomorphisms
$\a \ltimes_\al \f^\perp \cong \D(\g) \cong \g + \g^\st$ is given by
$d \mt (x(d), \x(d)).$ This implies that the bilinear form \eqref{in} on
$\a \ltimes_\al \f^\perp \cong \D(\g)$ is given by
\[
\lg \lg d_1, d_2 \rg \rg := \lg \x(d_1), x(d_2) \rg
+ \lg \x(d_2), x(d_1) \rg, \, d_i \in \a \ltimes_\al \f^\perp.
\]
\ere

\sectionnew{Applications to Poisson--Lie groups}
Let $G$ be a finite dimensional connected, simply connected
quasitriangular Poisson--Lie group with
Lie bialgebra $\g.$ (The ground field $k$ is assumed to be   
$\Rset$ or $\Cset.)$ As an application of 
Theorems~\ref{textension} and \ref{tg^st}, we
explicitly construct the double and dual Poisson--Lie groups 
$\D(G),$ $G^\st$ of $G.$
(Recall from the introduction that they are the simply
connected Poisson--Lie groups with tangent Lie bialgebras 
$\D(\g)$ and $\g^\st,$ respectively.) We also describe the 
image of $G^*$ under the induced map $j^\st \colon G^\st \ra \D(G)$ from
the embedding of
$\g^\st$ into $\D(\g)$.

We start with some general remarks on group extensions.
Assume that $H$ is a Lie group and $V$ is an $H$-module. A normalized
analytic 2-cocycle
$\mcA \in Z^2(H, V)$ gives rise to the extension $H \ltimes_{\mcA} V$
which as a manifold is $H \times V$ and has multiplication law
\[
(h_1, v_1) (h_2, v_2)= (h_1 h_2, v_1 + h_1 \cdot v_2 + \mcA(h_1, h_2) ).
\]
If $\h = \Lie(H)$ and the derivative of $\mcA$ is $\al \in Z^2(\h, V),$  
then $\Lie(H \ltimes_{\mcA} V) \cong \h \ltimes_{\al} V$. When $H$ is
simply connected any such Lie algebra cocycle $\alpha$ can
be lifted (uniquely) to an analytic Lie group cocycle $\mcA$ \cite{Ho}.

Abstract subgroups of $H \ltimes_{\mcA} V$ can be described in a
fashion analogous to \thref{subalgebras}. Let $B$ be a subgroup of
$H$ and let $\tW$ be a subgroup of $V$ that is $B$-invariant.
Denote by $\widehat{\sig}$ the projection $V \to V/\tW$.
A 1-cochain $\mcB\colon B \to V/\tW$ with coboundary
$- \widehat{\sig} \circ \mcA|_B$ gives rise to the subgroup
\[
 B^\mcB_\tW = \{ (b, v) \mid b \in B, \; v+\tW = \mcB(b) \}
\]
of $H \ltimes_{\mcA} V.$
Analogously to \thref{subalgebras}, all
abstract subgroups of $H \ltimes_{\mcA} V$ are obtained in this way. If
$B$ is a Lie subgroup, $\tW$ is closed in $V$ and $\mcB$ is analytic,
then $B^\mcB_\tW$ is a Lie subgroup.
 Let $\b$ be the Lie algebra of $B$ and let $W$ be the connected
component of $\tW$ (which is of course a subspace of $V$). The derivative
of $\mcB$ is a map $\be \colon \b \to V/W$ with coboundary
$- \sig \circ \al|_\b$.
In the notation of \leref{b-beta-W}, $\Lie(B^\mcB_\tW) \cong \b^\be_{W}$.

The converse is a little more subtle. Let $\b$ be a Lie subalgebra of $\h$
and let $B$ be the connected subgroup of $H$ corresponding to it.
Fix a 1-cochain $\be$ of $\b$ with values in $V/W$ and
coboundary $- \sig \circ \al|_\b$ and let $\b^\be_{W}$ be the
corresponding subalgebra of $\h \ltimes_{\al} V$. Denote by
$\eta \colon \tB \to B$ the universal cover of $B$, and let $\mcA_{\tB}$
be the induced 2-cocycle on $\tB$.
\ble{subgroups} 
The cochain $\beta$ of $\b$ can be integrated
to a unique 1-cochain $\tmcB$ on $\tB$ with values in $V/W$ and
coboundary $-\sig \circ \mcA_{\tB}.$ The subgroup
$\tB^\tmcB_{W}$ of $\tB\ltimes_{\mcA_{\tB}}V$ is then a
simply connected Lie group with Lie algebra $\b^\be_{W}$. The embedding
$\b^\be_{W} \hookrightarrow \h \ltimes_{\al} V$ yields a map 
$j^\st \colon \tB^\tmcB_{W} \to H \ltimes _\mcA V$ given by
$j^\st( \wt{b} \ltimes_{\mcA_{\tB}} v) = 
\eta(\wt{b}) \ltimes_\mcA v.$ Its image is the connected
subgroup of $H \ltimes_\mcA V$ with Lie algebra $\b^\be_{W}$ and is 
of the form
$B^\mcB_{\tW}$ where $\tW = \sig^{-1}(\mcB(\Ker \eta))$.
\ele

 
\noindent
\proof Denote by $K$ the connected subgroup of $\tB\ltimes_{\mcA_{\tB}}V$
with Lie algebra $\b^\be_{W}.$ Clearly $K \sup W$ and the projection
on the first component $p \colon K \ra \tB$ factors through to a 
projection $p_W \colon K/W \ra \tB.$ The induced map on Lie algebras
is an isomorphism. Since $\tB$ is simply connected, $p_W$ is also an
isomorphism and $K \cap V =W.$ Hence $K$ is of the type $\tB^\tmcB_{W}$
for a suitable map $\tmcB: \tB \ra V/W.$ Now $\tmcB$ must be an
analytic cochain on $\tB$ with coboundary $-\sig \circ \mcA_{\tB}$
because $K$ is a Lie subgroup of $\tB\ltimes_{\mcA_{\tB}}V.$
The uniqueness of a cocycle $\tmcB$ with the prescribed properties
follows from the uniqueness of a connected Lie subgroup with
given tangent Lie bialgebra.

The fact that $\tB^\tmcB_{W}$ is simply connected trivially follows from
the simply connectedness of $\Ker \, p_W=W$ and $\Im \, p_W= \tB.$
\hfill 
\qed 

Denote by $F$ the connected subgroup of $G$ with
$\Lie(F)= \f$ and by $G \ltimes F$ the semidirect product where
$G$ acts on $F$ by the adjoint action. Since $\f$ is an ideal, $F$ is
closed, normal and simply connected \cite[3.18]{Va}. Hence  $G\ltimes F$
is also simply connected. As remarked above the 2-cocycle $\al$ of $\g
\ltimes \f$ (see \eqref{bedef}) can be integrated to a 2-cocycle
$\mcA(., .)$ for the action of $G \ltimes F$ on $\f^\perp$ given by
\beq
(g, f) \cdot (\e) := \Ad^\st_g(\e)
\label{GFaction}
\eeq
The derivative of this action is of course the action of $\g \ltimes \f$ 
on $\f^\perp$ from \leref{ext-action}.

Consider the group
\[
D=(G \ltimes F) \ltimes_\mcA \f^\perp.
\]
Since
\[
\Lie(D) \cong (\g \ltimes \f) \ltimes_\al \f^\perp
\]
is isomorphic to the underlying Lie algebra of the double Lie bialgebra
$\D(\g),$ $D$
can be equipped with the Sklyanin Poisson structure (see, for instance,
\cite{CP}). The resulting Poisson--Lie group is the double of $G.$
As in the Lie algebra case, we can equivalently use the group
\[
A = \{ (g_1, g_2) \in G \times G \mid g^{-1}_1 g_2 \in F \}
\cong G \ltimes F
\]
and write $\mcA$ as a 2-cocycle for the action of $A$ on $\f^\perp$
induced from \eqref{GFaction} (cf. also \eqref{rep_d}).
Then
\[
D \cong A \ltimes_\mcA \f^\perp.
\]  

Denote by $G_\pm$ the connected subgroups of $G$ with tangent
Lie algebras $\g_\pm.$ Let $B$ be the connected subgroup
of $G_+ \times G_-$ with $\Lie(B) = \b.$ If
the isomorphism \\ $\th \colon \g_+/\n_+ \to \g_-/\n_-$
can be integrated to a group isomorphism $\Th \colon G_+/N_+ \to G_-/N_-$
for some closed (not necessarily connected) subgroups $N_\pm$ of $G$ with
$\Lie(N_\pm) = \n_\pm$
(taken the minimal possible), then $B$ is explicitly given by
\beq
B= \{ (g_1, g_2) \in G_+ \times G_- \mid \Th \circ \Pi_+ (g_1) =
\Pi_-(g_2) \}
\label{B}
\eeq
where $\Pi_\pm$ denote the homomorphisms $G_\pm \to G_\pm / N_\pm.$
As above $\eta \colon \tB \ra B$ will denote the simply connected cover 
of $B.$

Combining \thref{g^st} and \leref{subgroups} we obtain:

\bth{double} Let $G$ be a connected, simply connected
quasitriangular Poisson--Lie group. Let $F$ be the connected
normal Lie subgroup with Lie algebra $\f$. Then the double
Poisson--Lie group of $G$ is isomorphic, as a Lie
group, to $(G \ltimes F) \ltimes_\mcA \f^\perp \cong
A \ltimes_\mcA \f^\perp$ where $\mcA$ is a normalized analytic
2-cocycle on $G \ltimes F$ with derivative $\alpha$. The dual Poisson-Lie group 
$G^*$ is of the form  
$\tB_{W}^{\tmcB}$ where $W = \Ker \, r_+ \cap \Ker \, r_-$ and
$\tmcB: \tB \to V/W$ is the unique analytic 1-cochain on 
$\tB$ with coboundary $- \sig \circ \mcA_{\tB}$ and derivative
$\beta$.
The image of $G^*$ in $\D(G)$ is of the form $B^\mcB_{\tW}$ where
$\tW= \sig^{-1}(\tmcB(\Ker \eta))$ is a Lie subgroup of $V$ with
connected component $W$ and $\mcB: B \to V/\tW$ is a 1-cochain on $B$
induced from $\tmcB.$
\eth
\sectionnew{Special cases}
In this section we treat some cases of the results from the
previous ones and compare them to the available results in the literature.
\subsection{Factorizable Lie bialgebras}
Recall that a quasitriangular Lie bialgebra $(\g, r)$ is called
factorizable if $\Om = r + r_{21} \in S^2(\g)$ is nondegenerate,
or in other words if $\f=\g.$ In this case
$\f^\perp =0$ and $p$ defines an isomorphism
\[
\D(\g) \cong \a = \g \op \g.
\]
The results from Section~5 reduce to the fact that
in the above identification, $\g^\st \sub \D(\g)$ is isomorphic
to the subalgebra $\b \sub \g_+ \op \g_-,$ defined in \eqref{def_a}.
This recovers the following result
of Reshetikhin and Semenov-Tian-Shansky \cite{RS, STS}.

\bco{case1} The double of a factorizable Lie bialgebra
$\g$ is isomorphic, as a Lie algebra, to $\g \op \g.$
Under this isomorphism, the dual $\g^\st$ coincides with
the subalgebra $\b$ of $\g \op \g$ defined in \eqref{def_a} and $\g$ with
its diagonal subalgebra.
\eco

In particular \coref{case1} and Section~6 imply that the double of any
simply connected factorizable Poisson--Lie group is isomorphic to 
$G \times G$ as a Lie group. The image of the dual Poisson--Lie group
of $G$ inside $\D(G)$ coincides with the group $B$ constructed
in Section~6, see \eqref{B}.
\subsection{Triangular Lie bialgebras}
Let $(\g, r)$ be a triangular Lie bialgebra, that is
$r_{21}=-r.$ This means that $\Om=0$ and consequently
$\f=0,$ $\f^\perp=\g^\st.$ Hence the map $\s$ is trivial
and the exact sequence
\eqref{exact} naturally splits $(\al=0).$

\bco{tr1} The double $\D(\g)$ of a triangular Lie bialgebra $\g$
is isomorphic as a Lie algebra to the semidirect
sum of $\g$ and an abelian Lie algebra $\g^\st$
\beq
\D(\g) \cong \g \ltimes \g^\st
\label{triang_doub}
\eeq
where $\g$ acts on $\g^\st$ via the coadjoint action.
\eco

Next we specialize the results from Section~5 to the triangular
case. The Lie algebras coincides with $\g_+=\g_-.$
The space $W$ is $\Ker \, r_+$ and thus $\g^\st / W \cong \g^\st_+.$
The map $r_+ \colon \g^\st \to \g_+$ factors through to an
isomorphism $\ol{r}_+ \colon \g^\st/W \to \g_+.$ Its inverse
\beq
\be:= \ol{r}_+^{-1} \colon \g_+ \to \g^\st / \Ker \, r_+ \cong \g^\st_+ 
\label{cocyc_triang}
\eeq
is a 1-cocycle for the coadjoint representation of $\g_+.$

\bco{tr2} For a triangular Lie bialgebra $\g,$ the image of the
embedding $\g^\st \sub \D(\g) \cong \g \ltimes \g^\st$
coincides with
\[
(\g_+)^\be_{\Ker \, r_+} = (\id, \sig^{-1} \circ \be)(\g_+) \sub 
\g_+ \ltimes \g^\st.
\]
where $\be$ is the 1-cocycle \eqref{cocyc_triang} for the
coadjoint representation of $\g_+.$
\eco
As in Section~5, 
$\sig \colon \g^\st \to \g_+^\st \cong \g^\st/\Ker \, r_+$ denotes
the natural projection.
Notice that the image of  $\g^\st$ in $\D(\g)$ can be more simply described as
\[
\{(r_+(\xi), \xi) \mid \xi \in \g^\st \} \sub \g \ltimes \g^\st.
\]
However, this description, although explicit, is less useful because 
it cannot be lifted directly to the level of Poisson--Lie groups.

 Let $G$ be a simply connected triangular Poisson--Lie group and
$\g=\Lie(G)$. Then the double Poisson--Lie group of $G$
is isomorphic, as a Lie group, to $G\ltimes \g^\st.$ 
The dual Poisson--Lie group $G^\st$ of $G$ is isomorphic to
$(\widetilde{G}_+)^\tmcB_{\Ker \, r_+}$ where 
$\eta: \widetilde{G}_+ \ra G_+$ is the simply connected cover
of $G_+$ and $\tmcB \colon \widetilde{G}_+ \ra \g^\st/\Ker r_+$ is the
unique
1-cocycle on $\widetilde{G}_+$ with
derivative $\beta.$ The image of $G^\st$ in $G\ltimes \g^\st$ will be of
the form $(G_+)^\mcB_\tW$ where $\tW=\tmcB(\Ker \eta)+ \Ker r_+$ and 
$\mcB: G_+ \ra \g^\st/\tW$ is the 1-cochain on $G_+$ induced from $\tmcB.$ 
These results can be easily restated in terms of
$T^\st G \cong G \ltimes \g^\st.$

Notice that if $\beta$ is a coboundary 
(i.e., of the form $\beta(x)=ad^\st_x(\xi)$ for some 
$\xi \in \g_+^\st)$ then 
$\tmcB \colon \wt{G}_+ \ra \g_+^\st \cong \g^\st /\Ker r_+$ is
just the map $\tmcB(g) =(\Ad^\st_g-\id)\xi.$ 
Since $\Ker \/ \eta$ acts trivially on $\g^\st,$ $\tW=\Ker \/ r_+$ 
and $\mcB \colon \wt{G}_+ \ra \g_+^\st \cong \g^\st /\Ker r_+$ is given by
the same formula as $\tmcB$. 

Given the (nondegenerate) 1-cocycle $\be$ of $\g_+$ with values in
$\g_+^\st,$ one can define a
(nondegenerate) 2-cocycle $\ga(., .)$ of $\g_+$ with values
in the ground field $k$ defined by
\[
\ga(x_1, x_2) = \lg \be(x_1), x_2 \rg, \, x_i \in \g_+.
\]
This makes the pair $(\g_+, \ga)$ a quasi-Frobenius Lie algebra,
see \cite{Sto}. When $\be$ is nondegenerate, so is $\ga$ and the pair $(\g_+, \ga)$ becomes a Frobenius Lie algebra.

In the special case of a triangular bialgebra structure on a real compact
form $\g$ of a complex simple Lie algebra $\g_\Cset,$
Corollaries \ref{ctr1} and \ref{ctr2} were obtained by
Levendorskii and  Soibelman \cite{LS}. Although in this case
additional simplifications occur $(r \in \h \wedge \h$
for a Cartan subalgebra $\h$ of $\g),$ their proof can be extended
to the general case of triangular bialgebras. As
the one by Reshetikhin and Semenov-Tian-Shansky \cite{RS},
it relies on a direct construction of a Manin triple. \\
{\em{Sketch of a second proof of \coref{tr1}.}} Let $\d$ denote
the semidirect sum $\g \ltimes \g^\st,$ defined in \coref{tr1}.
Equip it with the following natural invariant bilinear form
\[
\lg (x_1, \x_1), (x_2, \x_2) \rg
= \lg \x_1 , x_2 \rg +  \lg \x_2 , x_1 \rg, \;
\forall x_i \in \g, \, \x_i \in \g^\st.
\]
It is clear that this form is nondegenerate and $\ad$ invariant.
Define the linear maps:
\beqa
&& j \colon \g \to \d,      \, j(x)=(x, 0),
\nn 
\\
&& j^\st \colon \g^\st \to \d,  \, j^\st(\x)=(r_+(\x), \x).
\nn
\eeqa
By a direct computation one checks that $j$ and $j^\st$ are
homomorphisms of Lie algebras and the images $j(\g) \cong \g,$
$j^\st(\g^\st) \cong \g^\st$ are Lagrangian subalgebras of
$\d.$ Thus $(\d, j(\g), j^\st(\g^\st))$ is a Manin triple
and $\D(\g) \cong \d$ as a Lie algebra.
\hfill \qed

If the Lie algebra $\g$ admits an invariant bilinear form,
$\g^\st$ can be identified, as a vector space, with $\g$ and the coadjoint
action substituted with the adjoint one. In this case
\coref{tr1} simply means that
\beq
\D(\g) \cong \g[t]/\lg(t^2)\rg
\label{triang_doub2}
\eeq
where $t$ is a free variable.
In addition $W= \Ker \, r_+$ can be identified with  $\g_+^\perp$
and $\sig$ thought of as projection $\g \to \g/ \g_+^\perp.$
Then $\be$ is a 1-cocycle of $\g_+$ with values in 
$\g / \g_+^\perp$ made a $\g_+$ module using the adjoint
action. (Everywhere the orthogonal complement is taken with 
respect to the nondegenerate form on $\g.)$ 
\coref{tr2} can be rephrased as follows.

\bco{trbf}If
$\g$ is a triangular Lie bialgebra admitting a nondegenerate
bilinear form, then
\[
\D(\g) \cong \g[t]/ \lg (t^2) \rg
\]
and the image of the embedding $\g^\st \sub \D(\g)$ is
\[
(\id + t \sig^{-1} \circ \be)(\g_+) \sub
\g[t]/ \lg  (t^2) \rg.
\]
\eco

Any  complex simple triangular Lie algebra
is of the type considered in \coref{trbf}.
In this case we recover a result
of Stolin, \cite{Sto}. The proof in \cite{Sto}
relies on works of Berman, Moody and
Benkart, Zelmanov on Lie algebras graded by
finite root systems. One of the motivations for this work was to
show that Corollaries \ref{ctr1} and \ref{ctr2}
follow from the triangularity of the bialgebra
structure without assuming simplicity of the algebra.
    
\end{document}